# TRIGONOMETRIC SPLINES IN SPECTRAL ANALYSIS PROBLEMS

Denysiuk V.P.
Dr. of Phys-Math. Sciences, Professor, Kiev, Ukraine
National Aviation University
kvomden@nau.edu.ua

Annotation
Some questions of application of trigonometric splines in problems of spectral analysis are considered. The known effects of overlay in the frequency and time domains are discussed; deployment effects in these areas are firstly considered. The relation between discrete Fourier coefficients and Fourier coefficients of trigonometric splines was obtained. The expediency of taking into account the differential properties of the investigated signals that are losing during the primary sampling is shown.

## Introduction

Currently, there are a lot of methods for spectral analysis of signals that are based on one or another system of assumptions about the investigated signal [1]. In this article we restrict ourselves to the classical spectral analysis of signals, which is based on the Fourier transform in space [2]. In particular, it is known that since the analog signal is always observed on a finite period of time, its spectrum is infinite and is represented by the Fourier coefficients [3].

The exact analog signal is either unknown or complex. Therefore, it is discretized, that is, it is replaced by a finite sequence of its discrete readings on a uniform grid.

The sampling process is accompanied by the well-known effect of overlay in the frequency (time) regions, for which the information on signal smoothness is required to evaluate the effect

After sampling most often consider the following approaches:

a) an approach in which the discrete Fourier coefficients are calculated for discrete countings by quadrature methods; the number of these coefficients are equal.

b) approach based on the orthogonality of trigonometric functions on the set of discrete samples; the number of discrete Fourier coefficients equal . This approach is well developed and covered in the literature [4], [5], [6].

c) approach based on the construction of the interpolation trigonometric polynomial, and get the discrete Fourier coefficients.

All three approaches yield similar results, which are the basis periodogram analysis [1].

There is a fundamentally different approach, based on calculation of integrals by the method of Philo. This approach is that in discrete counts of signal build some smooth function which approximates the signal in either sense; then compute the Fourier coefficients of this function.

Since in this approach it is possible to calculate countable many of Fourier coefficients, we can assume that in this case, the effect of deployment frequencies, which, generally speaking, is the opposite effect of imposing frequencies. It is clear that when constructing such an approximating function, it is necessary to take into account the information on the smoothness of the signal, which determine the order of decreasing coefficients and which are lost during sampling. In addition, it is desirable to provide the best approximation of the signal by its discrete countings.

This approach was proposed in [7]; however, no relation of the obtained Fourier coefficients with discrete coefficients was investigated.

It is possible to take into account the smoothness properties of the signal and at the same time ensure its best approximation using polynomial splines. A certain disadvantage of such splines is the complexity of constructing high-degree splines and their lumpy structure, which complicates the calculation of the Fourier coefficients of these splines.

In [8] it was built by classes of trigonometric splines, which involve polynomial's class. Trigonometric splines are given evenly by convergent Fourier series, which removes the need to calculate their Fourier coefficients. In addition, the connection between discrete Fourier coefficients and the Fourier coefficients of trigonometric splines is made, in which the differential properties of the investigated signal appears explicitly.

Special attention attracts the question of comparison of discrete Fourier coefficients and the Fourier coefficients of trigonometric splines. A direct comparison of these coefficients is unlikely to be appropriate, since each of the factors carries information about the signal that applies to the whole region of definition of the signal. Considering these factors, it is impossible to determine where the signal is, for example, a maximum or break (2). We have proposed to compare these factors as a whole in the time domain, that is, compare the quality of signal playback in the space metric $L^2$.

## The purpose of the work

Investigation of some issues related to the use of trigonometric splines in problems of spectral analysis.

## The main part

Class $2\pi$ - periodic functions with absolutely continuous derivative order $r-1$ (C $r=1,2,...$), and $r$ is the derivative of which is a function of bounded variation, let's denote it by $W_v^r$. Then, we shall call the functions of this class as signals

Signal $f(t) \in W_v^r$ you can submit it in a Fourier look

$$f(t) = \frac{a_0}{2} + \sum_{k=1}^{\infty} a_k \cos kt + b_k \sin kt \tag{1}$$

and the coefficients Фур'є $a_0$, $a_k$, $b_k$, calculated by the formulas $a_k = \frac{1}{\pi} \int_0^{2\pi} f(t)dt$;

$$a_k = \frac{1}{\pi} \int_0^{2\pi} f(t)\cos kt\, dt; \qquad b_k = \frac{1}{\pi} \int_0^{2\pi} f(t)\sin kt\, dt; \quad (k=1,2,...). \tag{2}$$

have a decreasing order $O(k^{-(1+r)})$. In addition, the coefficients have the following estimates

$$|a_k|, |b_k| \leq \frac{1}{\pi} \frac{\mathrm{Var}_0^{2\pi} f^{(r)}}{k^{r+1}}, \tag{3}$$

where $\mathrm{Var}_0^{2\pi} f^{(r)}$ - is a variation of the derivative $r$-th order of the signal $f(t)$.

When approximating the Fourier coefficients the signal is discretize, that is, replace the sequence of values at the nodes of the grid $\Delta = \{t_j\}_{j=1}^N$, $t_j = \frac{2\pi}{N}(j-1)$, $N=2n+1$, ($n=1,2,...$). It is considered that if the sample is sampled at a sufficiently high frequency, then the signal can be restored to any degree of accuracy. Below we show that such a statement is false.

It is clear that on $N$ discrete values $f(t_j)$ of signal $f(t)$ in nodes of grid $\Delta$ we can get only $N$ independ discrete coefficients $a_0^*$, $a_k^*$, $b_k^*$ ($k=1,2,...n$), which are calculated by formula

$$a_0^* = \frac{2}{N}\sum_{j=1}^N f(t_j)$$

$$a_k^* = \frac{2}{N}\sum_{j=1}^N f(t_j)\cos kt_j; \qquad b_k^* = \frac{2}{N}\sum_{j=1}^N f(t_j)\sin kt_j, \quad (k=1,2,...,n). \tag{4}$$

These formulas are easy to obtain by quadrature methods or by assuming the orthogonality of trigonometric functions on a discrete set of points, or by calculating the coefficients of an interpolation trigonometric polynomial. Note that although the discrete coefficients $a_k^*$, $b_k^*$ can be calculated for any index value $k$, $1 < k < \infty$, thus receive only $N$ the periodic, even, and odd sequences of these coefficients.

When sampling the signal $f(t)$ there is a well-known frequency-mapping effect.

We calculate the value of the signal $f(t)$ in $j+1$-th node of a grid $\Delta$; because of (1) consider:

$$f(t_{j+1}) = \frac{a_0}{2} + \sum_{k+1}^{\infty}(a_k \cos kt_{j+1} + b_k \sin kt_{j+1}) = \frac{a_0}{2} + \sum_{k=1}^{\infty}\left(a_k \cos k\frac{2\pi}{N}j + b_k \sin k\frac{2\pi}{N}j\right) \tag{5}$$

Because the value of the functions $\cos\frac{2\pi}{N}kj$ and $\cos\frac{2\pi}{N}(mN\pm k)j$ (the same for sinuses), in nodes of grid $\Delta$ coincide, (5) can be submitted as

$$f(t_{j+1}) = \left[\frac{a_0}{2}+\sum_{m=1}^{\infty}a_{2mN}\right] + \sum_{k=1}^{n}\cos k\tfrac{2\pi}{N}j\left[a_k + \sum_{m=1}^{\infty}(a_{mN+k}+a_{mN-k})\right] +$$
$$+ \sum_{k=1}^{n}\sin k\tfrac{2\pi}{N}j\left[b_k + \sum_{m=1}^{\infty}(b_{mN+k}-b_{mN-k})\right] \qquad (6)$$

Introducing the notation

$$\frac{a_0^*}{2} = \frac{a_0}{2} + \sum_{m=1}^{\infty}a_{2mN};$$
$$a_k^* = a_k + \sum_{m=1}^{\infty}(a_{mN+k}+a_{mN-k}); \qquad (7)$$
$$b_k^* = b_k + \sum_{m=1}^{\infty}(b_{mN+k}-b_{mN-k}),$$

equation (6) looks like

$$f(t_{j+1}) = \frac{a_0^*}{2} + \sum_{k=1}^{n}a_k^*\cos kt_{j+1} + b_k^*\sin kt_{j+1},$$

Where coefficients $a_0^*$, $a_k^*$, $b_k^*$ ($k=1,2,...n$) are calculated by formulas (4).

Ratio (7) is the so-called frequency overlay effect (look example [4], [5], [6]), that occurs when a signal is replaced $f(t)$ the sequence of its values in the nodes of a uniform grid $\Delta_N$.
From (3) it follows that

$$|a_k - a_k^*| \leq \sum_{m=1}^{\infty}[|a_{mN+k}| + |a_{mN-k}|] \leq$$
$$\leq \frac{1}{\pi}\mathop{Var}_0^{2\pi}f^{(r)}\sum_{m=1}^{\infty}\left[\frac{1}{(mN+k)^{r+1}} + \frac{1}{(mN-k)^{r+1}}\right]. \qquad (8)$$

A similar estimate holds for the coefficients $b_k$ and $b_k^*$.

So far, we have considered the effect of overlay in the frequency domain. However, it is often convenient to consider this effect in the time domain [2].

Signal component $f_n(t)$, which is transported to $N$ as first Furier's coefficients $a_0$, $a_k$, $b_k$ ($k=1,2,...n$) served by expression

$$f_n(t) = \frac{a_0}{2} + \sum_{k=1}^{n}(a_k\cos kt + b_k\sin kt).$$

Accordingly, the component of the signal $f_n^*(t)$, which is transported by the higher Fourier coefficients, has the form (we do not take into account the constant component, which does not play a significant role in the problems of spectral analysis)

$$f_n^*(t) = \sum_{m=1}^{\infty}[a_{mN+k}\cos(mN+k)t + a_{mN-k}\cos(mN-k)t] +$$
$$+ \sum_{m=1}^{\infty}[b_{mN+k}\sin(mN+k)t + b_{mN-k}\sin(mN-k)t],$$

So the signal $f^*(t)$, , which is given by discrete Fourier coefficients and has the form

$$f^*(t) = \frac{a_0^*}{2} + \sum_{k=1}^{n}a_k^*\cos kt + b_k^*\sin kt$$

can be submitted as follows

$$f^*(t) = f_n(t) + f_n^*(t).$$

Given ratio and gives the effect of overlay in the time domain.

Because of the function $f_n^*(t)$ can be estimated in the standard way

$$|f_n^*(t)| \leq \frac{2}{n^r}\mathop{Var}_0^{2\pi}f^{(r)},$$

we finally have

$$\left|f_n(t)-f^*(t)\right|\leq \frac{2}{n^r}\underset{0}{\overset{2\pi}{Var}} f^{(r)}. \qquad (9)$$

From (8), (9) it is important to conclude that the estimates of the overlay effect in the frequency and time domains depend significantly on the order of decreasing Fourier coefficients; in turn, the order of decrease of these coefficients is determined by the differential properties of the investigated signal, які втрачаються when switching from analog signal to discrete; it is impossible to account for these properties under formulas (4). These differential properties can also be taken into account by calculating the approximate Fourier coefficients by the Filon method. Let's consider this method more detailed. The Philon method, proposed for calculating integrals of type (2) is as follows. Instead of approaching the entire integrand, only the signal approaches $f(t)$, but rather fast shaking factors $\cos kx$ and $\sin kx$ are considered as weight factors.

This approach to the approximate calculation of integrals of type (2) has several advantages over what we have previously considered. Among these advantages, we will recall only the following.

First, when calculating integrals of type (1) with the Filon method, we need to approximate only the signal $f(t)$, which changes relatively slowly; accordingly, some computational accuracy can be achieved with a much smaller number of interpolation grid nodes.

Secondly, the approximate Fourier coefficients calculated by the Filon method, are ordinary Fourier coefficients of the approximate function and can be investigated by known methods.

Finally, the main advantage of the Filon method is that we can calculate Fourier coefficients with any indexes. In other words, there is a frequency deployment effect that is opposite to the frequency overlay effect.

Consider an approximate Fourier coefficient calculation by the Filon method.

Let $\varphi_N(t)$ - some function approaching the signal $f(t)$ in one sense or another, using a sequence of discrete values of this function at the nodes of the grid $\Delta$. Then the approximate Furier coefficients are calculated by the formulas $\hat{a}_k = \frac{1}{\pi}\int_0^{2\pi}\varphi_N(t)dt$;

$$\hat{a}_k = \frac{1}{\pi}\int_0^{2\pi}\varphi_N(t)\cos kt\, dt; \quad \hat{b}_k = \frac{1}{\pi}\int_0^{2\pi}\varphi_N(t)\sin kt\, dt; \quad (k=1,2,\ldots),$$

moreover, the number of these coefficients is, in general, unlimited. This, in fact, is the deployment effect which we mentioned earlier.

Consider the estimation of the error of the approximate Fourier coefficients calculated by the Filon method. This error can be estimated by the formula

$$\left|a_k - \hat{a}_k\right| \leq \max_t \frac{1}{\pi}\int_0^{2\pi}\left|f(t)-\varphi_N(t)\right|\left|\cos kt\right| dt = \frac{1}{\pi}\left\|f(t)-\varphi_N(t)\right\|_C \int_0^{2\pi}\left|\cos kt\right| dt =$$

$$= \frac{1}{\pi}\left\|f(t)-\varphi_N(t)\right\|_C 4k\int_0^{\frac{\pi}{2k}}\cos kt\, dt = \frac{4}{\pi}\left\|f(t)-\varphi_N(t)\right\|_C.$$

Thinking similarly, we get:

$$\left|b_k - b_k^*\right| \leq \frac{4}{\pi}\left\|f(t)-f^*(t)\right\|_C.$$

Thus, the error in the calculation of the Fourier coefficients using the Filon method depends only on the error of the approximation of the signal $f(t)$. As it is known [9 that the best linear apparatus of approximation (whose classes are functions) are simple polynomial (and, accordingly, trigonometric) splines, we can conclude that in order to ensure the least error in the calculation of the Fourier coefficients in the calculation of these coefficients, it is advisable to use these splines. Below we obtain a similar conclusion based on other considerations.

The disadvantage of this estimate is that it does not depend on the value of the parameter. Because with the growth $k$ the Fourier coefficients of continuous functions decrease, the obtained estimates do not guarantee a small relative error of the higher Fourier coefficients.

This situation can be corrected as follows. Let's return to the consideration of quantities

$$a_k - \hat{a}_k = \frac{1}{\pi}\int_0^{2\pi}\left[f(t)-\varphi_N(t)\right]\cos kt\, dt$$

$$b_k - \hat{b}_k = \frac{1}{\pi}\int_0^{2\pi}\left[f(t)-\varphi_N(t)\right]\sin kt\, dt.$$

It is clear that the integrals on the right-hand side of the obtained equations are Fourier coefficients of the difference between the two functions and can be investigated by standard methods.

For example, let's have a function $f(t) \in W_\nu^r$ and approaching function $\varphi_N(t) \in W_\nu^m$. It is clear that under these assumptions there is a difference $[f(t) - \varphi_N(t)] \in W_\nu^q$, ($q = \min(r,m)$), and accordingly, there is an evaluation $|a_k - \hat{a}_k| \leq \dfrac{1}{\pi k^{q+1}} \operatorname*{Var}_0^{2\pi}[f(t) - \varphi_N(t)]^{(q)}(t)$,

$$|b_k - \hat{b}_k| \leq \dfrac{1}{\pi k^{q+1}} \operatorname*{Var}_0^{2\pi}[f(t) - \varphi_N(t)]^{(q)}(t).$$

The following approach to estimating the Fourier coefficient's calculation error makes the following conclusion.

The order of fall of the Fourier coefficient calculation error is determined by the differential properties as a function $f(t)$ the same as function $f^*(t)$, and can't be higher then $O\left(\dfrac{1}{n^{q+1}}\right)$, where $q = \min(r,m)$. calculation error is achieved when the differential properties of the function $\varphi_N(t)$ are not lower than the differential properties of the signal $f(t)$.

The conclusion has rather qualitative character, because it is possible to obtain estimates of the variation of derivative differences $[f(t) - \varphi_N(t)]$ is quiet difficult. However, this conclusion can be used to select a class of approximate functions.

Construct an approximate function $\varphi_N(t)$, which has certain differential properties and provides the best approximation of the signal using trigonometric splines [8]. All the results obtained for these splines are also extended in the case of their corresponding polynomial periodic splines.

Let us consider the Philo method of calculating approximate Fourier coefficients for the case when interpolating simple trigonometric splines are used as approximate functions.

In [8] it was shown that trigonometric spline $St(r,t)$ of an order $r$, ($r = 1, 2, \ldots$), which interpolates the signal at the grid nodes $\Delta$, can be submitted as

$$St(r,t) = \dfrac{a_0^*}{2} + \sum_{k=1}^{N-1} \dfrac{1}{H(r,k)} \sum_{m=0}^{\infty} \sigma_{mN+k}(r)\left[a_{mN+k}^* \cos(mN+k)t + b_{mN+k}^* \sin(mN+k)t\right]$$

where

$a_0^*$, $a_{mN+k}^*$, $b_{mN+k}^*$, - discrete Fourier coefficients calculated by formulas (4) for the respective indices;

$\sigma_{mN+k}(r)$ - factors of type $\varphi(r,k)/k^{1+r}$, in descending order $O(k^{-(1+r)})$, and $\varphi(r,k)$ is some limited function, and a factor $H(r,k)$ is calculated by the formula

$$H(r,k) = \sigma_k(r) + \sum_{m=1}^{\infty}\left[\sigma_{mn+k}(r) + \sigma_{mn-k}(r)\right].$$

Note, that in the role of factors $\sigma_j(r)$ you can choose, for example, factors

$$\sigma_j(r) = \left[\operatorname{sinc}(\dfrac{\pi}{N} j)\right]^{1+r}, \quad \sigma_j(r) = \left|\operatorname{sinc}(\dfrac{\pi}{N} j)\right|^{1+r}, \quad \sigma_j(r) = \dfrac{1}{j^{1+r}},$$

$$(j = 1, 2, \ldots).$$

where checked

$$\operatorname{sinc}(x) = \sin x / x.$$

It is clear that the trigonometric spline $St(r,t) \in W_\nu^r$, and its Furier coefficients have the form

$$\hat{a}_j = \dfrac{\sigma(r,j)}{H(r)} a_j^*; \qquad \hat{b}_j = \dfrac{\sigma(r,j)}{H(r)} b_j^*, \qquad (j = 1, 2, \ldots).$$

Therefore, in this case there is no need to calculate the Fourier coefficients since, as we have already said, the spline itself is supplied by the Fourier series.

Particular attention paid to the factors $\alpha(r,j) = \sigma(r,j)/H(r)$, which are included in the spline coefficients. Here are graphs of such factors for some parameter values $r$ for the case of factors $\sigma_j(r)$ and factor $\sigma_j(r) = \left[\operatorname{sinc}(\dfrac{\pi}{N} j)\right]^{1+r}$ was selected.

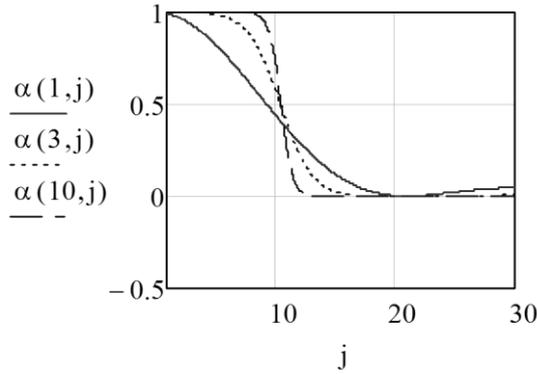

Picture 1. Graphs of quantities $\alpha(r, j)$ for values $r = 1, 3, 10$.

In the graphs above it is easy to find the graphs of the amplitude-frequency characteristics of the low-pass filters [10], the cutoff frequency of which is determined by the Nyquist frequency, and the slope's decay is determined by the differential properties of the signal. This, allows to substantiate the empirical methods of construction of low-pass filters, which are often used in practice. In addition, theoretical calculations open the possibility of taking into account frequencies which are higher then Nyquist frequencies.

Let us now turn to the analysis of the effect of deployment in the time domain. Since polynomial and, accordingly, trigonometric splines are supposed to be used as approximating functions, it is advisable to refer to their extreme properties. Yes, in particular, it is known that a derivative $2k$-th order of polynomial spline $2k+1$-th order has the property of the slightest curvature ($k = 0, 1, ...$). This property means that all functions which have a continuous derivative $2k$-th order and interpolates a given sequence of counts of the signal, only the spline derivative minimizes the functional

$$\int_0^{2\pi} \left[ f^{(2k)}(t) \right]^2 dt.$$

Since polynomial spline order $2k+1$ coincide with trigonometric splines of order $r$, ($r = 2k+1$), the minimum curvature property can also be transferred to trigonometric splines. From this follows that the variation of this derivative, which is included in the estimate (9), takes the smallest value.

## Conclusions

1. The role of differential properties that are lost during sampling in signal processing problems is considered.
2. Estimates of the harmful effect of the frequency overlay effect using the differential properties of signals in the frequency domain
3. An estimation of the harmful effect of the frequency overlay effect in the time domain is given.
4. The feasibility of using interpolation polynomial and trigonometric splines in the calculation of approximate Fourier coefficients by the Filon method.
5. The effect of frequency unfolding observed in the calculation of approximate Fourier coefficients by the Filon method is firstly considered.
6. The relation between the discrete Fourier coefficients and the approximate Fourier coefficients calculated by the Filon method is established.
7. The empirical methods of construction of low-pass filters, which are often applied in practice, are substantiated; in addition, theoretical calculations make it possible to consider frequencies then are higher then Nyquist frequencies.
8. The results of numerical simulation are in good agreement with the theoretical results.
9. 9. Of course, the obtained results require further investigation.

## List of references.